\documentclass[twoside,11pt]{amsart}
\usepackage{epsfig}
\usepackage{amsfonts}
\usepackage{amssymb}

\renewcommand{\baselinestretch}{1.2}

\begin{document}
\title{INTERPOLATION MAPS AND CONGRUENCE DOMAINS FOR WAVELET SETS}
\author{Xiaofei Zhang and David R. Larson}
\address{Department of Mathematics, Texas A$\&$M University, College
Station, TX 77843; xfzhang1974@yahoo.com; larson@math.tamu.edu}

\date{}

\
\thanks{*The second author was partially supported by NSF grant DMS-0139386.}

\begin{abstract}
It is proven that if an interpolation map between two wavelet sets
preserves the union of the sets, then the pair must be an
interpolation pair. We also construct an example of a pair of
wavelet sets for which the congruence domains of the associated
interpolation map and its inverse are equal, and yet the pair is not
an interpolation pair.  The first result solves affirmatively a
problem that the second author had posed several years ago, and the
second result solves an intriguing problem of D. Han.  The key to
this counterexample is a special technical lemma on constructing
wavelet sets.  Several other applications of this result are also
given.  In addition, some problems are posed.  We also take the
opportunity to give some general exposition on wavelet sets and
operator-theoretic interpolation of wavelets.

\end{abstract}

\maketitle

\begin{center}

{\it Dedicated to Larry Baggett for his great friendship, his love of
\\ mathematics, and his continued support of young mathematicians.}

\vspace{12pt}
\end{center}

\section{Introduction}
An orthonormal wavelet is a single function $\psi$ in
${\mathcal{L}}^2(\mathbb{R}^n)$ whose translates by all members of a
full-rank lattice followed by dilates by all integral powers of a
real expansive matrix on $\mathbb{R}^n$ generates an orthonormal
basis for ${\mathcal{L}}^2(\mathbb{R}^n)$. By the term \emph{wavelet
set} we mean a measurable subset $E \subset \mathbb{R}^n$ with the
property that the inverse Fourier Transform of its normalized
characteristic function, $\mathcal{F}^{-1}(\frac{1}{(2\pi)^{\frac{n}{2}}}\chi_E)$,
 is an orthonormal wavelet.
In [DL] an operator-theoretic technique for working with certain
problems concerning wavelets was introduced that was called {\it
operator-theoretic interpolation}. If $\psi$ and $\eta$ are
orthonormal wavelets in the same space, and $\psi_{n,l}$ and
$\eta_{n,l}$ are the corresponding wavelet bases, then the unitary
operator determined by the mapping $\psi_{n,l}$ to $\eta_{n,l}$
was called the {\it interpolation unitary} between $\psi$ and
$\eta$. These interpolation operators associated with ordered pairs
of wavelets play an essential role in the theory. They are
associated with the von Neumann subalgebras of the so-called
\emph{local commutant} space, whose unitary groups provide natural
parameterizations of certain families of wavelets. In the special
case where the interpolation operator is involutive (i.e. has square
I) the pair of wavelets is called an \emph{interpolation pair} of
wavelets. One surprising feature of the theory is that interpolation
pairs occur not infrequently.

In general, interpolation unitaries can be hard to work with.
However, in the special case where $\psi$ and $\eta$ are {\it
s-elementary wavelets} (also called {\emph MSF-wavelets} with phase
$0$, or wavelet-set wavelets), the interpolation unitary takes the
form of a \emph {composition operator} with measure-preserving
symbol. Every pair of wavelet sets gives rise to a measure
preserving transformation on the underlying measure space in a
natural way, called the {\it interpolation map} determined by the
pair. The interpolation theory for such wavelets and their
associated wavelet sets is, in many concrete cases, computable by
hands-on experimental paper-and-pencil computations. This permits
experimentation in the form of testing of hypotheses in potential
theorems for more general types of wavelets. The simplest case is
where a pair of wavelet sets has the property that the measure
preserving transformation $\sigma$ is an involution (i.e. $\sigma
\circ \sigma = id$). In this case the composition unitary has square
$I$, so the pair of wavelets is indeed an interpolation pair. The
pair of wavelets sets is, by analogy, called an {\it interpolation
pair of wavelet sets}. More generally, an interpolation family of
wavelets (and analogously, of wavelet sets) is a finite (or even
infinite) family of wavelets for which the associated family of
interpolation unitaries (interpolation maps) forms a group.

It is appropriate to give a bit of background and history that
serves to indicate why interpolation pairs of wavelet sets are
relevant to the theory of wavelets. More exposition on this,
including specific details and statements of theorems involved, can
be found in the semi-expository articles [La2], [La3], [La4], and
[La5]. In [DL], for any interpolation pair of wavelet sets (E, F),
the authors constructed a 2 x 2 complex matrix valued function
(called the Coefficient Criterion, see [DL], Proposition 5.4, and
also [La2], section 5.22) which specifies precisely when a function
$f$ on $\mathbb{R}$ with frequency support contained in $E \cup F$
is an orthonormal wavelet. If $supp\{\hat{f}\}$ is contained in
$E\cup F$, then this criterion shows that $f$ is an orthonormal
wavelet iff this matrix valued function is a unitary matrix (a.e.),
and it is a Riesz wavelet iff it is an invertible matrix (a.e).
Moreover, it shows that a Parseval frame wavelet (resp. Riesz frame
wavelet) with frequency support contained in $E \cup F$ is
necessarily an orthonormal wavlet (resp. Riesz wavelet). It follows
that the set of orthonormal wavelets with frequency support
contained in the union of $E \cup F$, where $(E,F)$ is an
interpolation pair of wavelet sets, is pathwise connected in
$\mathcal{L}^2(\mathbb{R})$. The set of Fourier Transforms of this
set is also connected in the $\mathcal{L}^{\infty}$ norm on the
frequency space.

These results were the main motivating factor in posing the first
open problem discussed in [DL], namely the question of whether the
set of all orthonormal wavelets in $\mathcal{L}^2(\mathbb{R})$ is
norm-pathwise connected. This was the same problem that was posed
completely independently by G. Weiss and his research group in
[HWW1],[HWW2] for different reasons. Their reasons included the
interesting discovery that certain wavelet sets (rather, the
associated MSF wavelets) could be "smoothed" in a continuous fashion
to obtain wavelets that were continuous in the frequency domain. It
turned out that our operator-interpolation approach, in certain key
cases, was equivalent to the "smoothing" approach of G. Weiss, and
the cases involved included the derivation of Y. Meyer's classic
family of wavelets that are compactly supported and continuous in
the frequency domain. Exploring common interest in the
relationships between smoothing of a wavelet set on the one hand,
and operator-theoretic interpolation between a pair of wavelet sets
on the other hand, and the general "connectedness" problem that was
motivated by both approaches independently (as described above) let
to the formation of the WUTAM Consortium (short for Washington
University and Texas A\&M University) and the joint work [Wut] of
the consortium, in which the connectedness problem was shown to have
a positive answer for the case of MRA wavelets.

The basic idea behind operator-interpolation is elementary. If $x$
and $y$ are elements of a vector space $V$, we say that a vector $z$
is linearly interpolated from $x$ and $y$ if $z$ is a convex
combination of $x, y$. More generally, it is convenient to allow
arbitrary linear combinations. So the set of vectors interpolated
by $x$ and $y$ is the linear span of $x$ and $y$. More generally,
we can say that $z$ is linearly interpolated from a collection
$\mathcal{F}$ of vectors if $z$ is a linear combination of vectors
from the family. And more generally yet, if the vector space $V$ is
a left module over some operator algebra $\mathcal{D}$ we can
consider linear combinations from $\mathcal{F}$ with coefficients
that are operators from $\mathcal{D}$, called \emph{modular} linear
combinations. If $z$ is a modular linear combination of $x, y$, then
we say that $z$ is derived from $x$ and $y$ by
operator-theoretic-interpolation (or operator-interpolation for
short). In the case of wavelets, the operator algebra $\mathcal{D}$
is the von Neumann algebra of all bounded linear operators acting on
$\mathcal{L}^2(\mathbb{R}^n)$ that commute with the dilation and
translation unitary operators for the wavelet system. When we
conjugate this with the Fourier transform (which is unitary), so we
are working in $\mathcal{L}^2(\mathbb{R}^n)$ as the frequency space,
and if we denote this conjugated algebra by $\hat{\mathcal{D}}$,
then $\hat{\mathcal{D}}$ is an algebra of \emph{multiplication
operators} on $\mathcal{L}^2(\mathbb{R}^n)$. In particular, it is a
commutative algebra. If a wavelet $\eta$ is a modular linear
combination of wavelets $\psi$ and $\nu$ with coefficients which are
operators in $\mathcal{D}$, then we say $\eta$ is derived by
operator-interpolation between $\psi$ and $\nu$. Not all modular
linear combinations of $\psi$ and $\nu$ are orthonormal wavelets.
They are all \emph{Bessel} wavelets to be sure, but a certain
unitarity condition needs to be satisfied to be an orthonormal
wavelet. Let ${V_{\psi}}^{\nu}$ be the interpolation unitary from
$\psi$ to $\nu$. If $A$ and $B$ are operators in $D$ and $\eta =
A\psi + B\nu$, the necessary and sufficient condition for $\eta$ to
be an orthonormal wavlet is that the operator $U := A +
B{V_{\psi}}^{\nu}$ needs to be unitary. (More generally, for a frame
wavelet the criterion is that $U$ must be surjective, and for a
Riesz wavelet $U$ must be invertible.)

The reason that interpolation pairs of orthonormal wavelets are
special is that if $(\psi, \nu)$ is an interpolation pair, and if
$A$ and $B$ are operators in $D$ such that $A*A + B*B = I$, then
under certain circumstances $A\psi + B\nu$ can also be an
orthonormal wavelet.
 In particular, if $\theta \in [0, 2\pi]$ is arbitrary, then
 $\eta := cos{\theta} \psi + i sin{\theta} \nu$
 is an orthonormal wavelet. Indeed, in this case the operator $U$ above
is just $cos{\theta} I + i sin{\theta} V_{\psi}^{\nu}$, and since
$({V_{\psi}}^{\nu})^2 = I $ it follows that $UU^* = I$, so $U$ is
unitary, as required by the criterion. Letting $\theta$ vary
continuously it follows, in particular, that $\psi$ and $\nu$ are
pathwise connected via a path of orthogonal wavelets.

For wavelet set wavelets, i.e. MSF wavelets with phase $0$, more is
true, and the operator-algebraic geometry involved is fairly rich.
For \emph{any} pair of wavelets sets $E$ and $F$, with associated
MSF wavelets $\psi_E$ and $\psi_F$, the interpolation unitary
$V_{\psi_E}^{\psi_F}$ \emph{normalizes} the von Neumann algebra
$\mathcal{D}$ in the sense that $V_{\psi_E}^{\psi_F} \mathcal{D}
(V_{\psi_E}^{\psi_F})^* = \mathcal{D}$. This was proven in Chapt 5
of [DL], and was a key result of that memoir. (We note that it is
an open question (see [La2], Problem 4) as to whether arbitrary
interpolation operators (i.e. for non-wavelet set wavelets)
normalize $\mathcal{D}$.) The reason that interpolation pairs of
wavelet sets, are even more special than general interpolations of
wavelets is the following: Firstly, $(E, F)$ is an interpolation
pair of wavelet sets if and only if the pair of wavelets $(\psi_E,
\psi_F)$ is an interpolation pair of wavelets. And secondly, since
in this case the interpolation unitary normalizes $\mathcal{D}$ ,
and since $({V_{\psi_E}}^{\psi_F})^2 = I$, it follows that the set
of operators $\{A + B V_{\psi_E}^{\psi_F} \, | \, A, B \in
\mathcal{D}\}$ is closed under multiplication and is in fact a von
Neumann algebra. (For a more general interpolation pair of wavelets
whose interpolation operator normalizes $\mathcal{D}$ the same thing
is true.) Since the unitary group of a von Neumann algebra is
pathwise connected in the operator norm, the interpolated family of
wavelets is also connected. Much more is true. Since the elements
of $\mathcal{D}$ are multiplication operators in a certain family
(the dilation-periodic operators), if we write $A = M_f$ and $B =
M_g$ we obtain $(A + B V_{\psi_E}^{\psi_F})\psi_E = f\psi_E +
g\psi_F$. Thus the outcome is an actual formula as well as a
criterion (the Coefficient Criterion mentioned above), for
constructing all orthonormal wavelets whose frequency support is
contained in the union $E \cup F$ of a given interpolation pair
$(E,F)$ of wavelet sets. By choosing $f$ and $g$ appropriately, so
they are continuous and vanish on the boundary of the union $E \cup
F$, and agree on $E \cap F$, and satisfy the unitarity condition
referred to above, one can obtain wavelets in this fashion which are
smooth in the frequency-domain. Not all interpolation pairs $(E,F)$
can be so \emph{smoothed}. But some can. As alluded to above, Y.
Meyer's famous class of orthonormal wavelets which are continuous
and compactly supported in the frequency domain can be derived in
this way from a special interpolation pair of wavelet sets: namely
the pair $E = [-\frac{8\pi}{3}, -\frac{4\pi}{3}) \cup
[\frac{2\pi}{3}, \frac{4\pi}{3})$ and $F = [-\frac{4\pi}{3},
-\frac{2\pi}{3}) \cup [\frac{4\pi}{3}, \frac{8\pi}{3})$.
 Even for cases in which
smoothing cannot work (for instance, $E \cup F$ may have too many
boundary points), the operator algebra involved can be interesting.

The purpose of this paper is to provide solutions to two related
problems concerning dyadic orthonormal wavelet sets in the line. One
problem asked whether a certain containment relation for an
interpolation map implies that the associated pair of wavelet sets
is an interpolation pair. This was posed by the second author in a
VIGRE seminar course at Texas A\&M several years ago. We answer
this question affirmatively, and we also observe that the analogous
result does not hold for a general interpolation family of wavelet
sets. The second problem was posed by D. Han, who asked whether the
equality of the congruence domains of an interpolation map and its
inverse implies that the associated pair of wavelet sets is an
interpolation pair. We were able to give a counterexample to this
problem. Our work on this interesting problem motivated a useful
lemma on constructing wavelet sets, which is apparently different
from those methods which have appeared in the literature to date,
and which is used in this counterexample as well as in the
construction of several wavelet sets concerning some related
questions, and also some wavelet sets in the plane. Much of the work
presented in this article is material from the doctoral dissertation
of the first author [Zh], which has not appeared elsewhere.

Much work has been accomplished on the topic of wavelet sets since
the mid 1990's. We have outlined some of the background and history
in the opening paragraph of Section 5, for the interested reader.
Our main results in this paper are for the special case n = 1 with
the dilation scale factor 2 and integer translates (the dyadic
case). In Section 5, we apply one of our techniques to the
construction of certain dyadic wavelet sets in the plane. Based on
our results, some further directions are suggested for higher
dimensions.

\section{Two Problems}

A {\it dyadic orthonormal wavelet} is a function $\psi \in
{\mathcal{L}}^2(\mathbb{R})$ (Lebesgue mesasure) with the property
that the set $\{ \, 2^{\frac{n}{2}}\psi (2^n\cdot -l) \, |\,
n,l\in\mathbb{Z}\, \}$ forms an orthonormal basis for
${\mathcal{L}}^2(\mathbb{R})$. More generally, if $A$ is any real
invertible $n\times n$ matrix, then a single function $\psi\in
L^2({\mathbb{R}}^n)$ is an orthonormal wavelet for $A$ if
\[ \{ |det A|^{\frac{n}{2}} \psi(A^n\cdot -l) \,|\, n\in \mathbb{Z}, l\in
\mathbb{Z}^{(n)} \} \] is an orthonormal basis of
$L^2({\mathbb{R}}^n)$. If $A$ is {\it expansive} (equivalently, all
eigenvalues of $A$ are required to have absolute value strictly
greater than $1$) then it was shown in [DLS1] that orthonormal
wavelets for $A$ always exist.

By the support of a measurable function we mean the set of points in
its domain at which it does not vanish. By the support of an element
$f$ of ${\mathcal{L}}^2(\mathbb{R})$, we mean the support of any
measurable representative of $f$, which is well-defined in the
measure algebra of equivalence classes of sets modulo null sets. By
the \emph{frequency support} of a function we mean the support of
its Fourier Transform.

Let $\mathcal{F}$ denote the $n$-dimensional Fourier transform on
$\mathcal{L}^2(\mathbb{R}^n)$ defined by
\[ (\mathcal{F} f)(s) := \frac{1}{(2\pi)^{\frac{n}{2}}}
\int_{\mathbb{R}^n} e^{-s\circ t} f(t) dm \]
for all $f\in L^2(\mathbb{R}^n)$. Here, $s\circ t$ denotes the real inner
product. A measurable set $E\subseteq \mathbb{R}^n$ is a wavelet set for
$A$ if
\[ \mathcal{F}^{-1} (\frac{1}{\sqrt{\mu(E)}} \chi_E) \]
is an orthonormal wavelet for $A$.

A sequence of measurable sets $\{ E_n \}$ is called a measurable
partition of $E$ if $E = E \,\Delta\, (\bigcup_n E_n)$ is a null set
and $E_n \cap E_m$ has measure zero if $n\ne m$, where $\Delta$
denotes the symmetric difference of sets. Measurable subsets $E$ and
$F$ of $\mathbb{R}$ are called {\it $2\pi$-translation congruent} to
each other, denoted by $E\sim_{2\pi} F$, if there exists a
measurable partition $\{ E_n\}$ of $E$, such that $\{ E_n + 2n\pi\}$
is a measurable paritition of $F$. Similarly, $E$ and $F$ are called
{\it 2-dilation congruent} to each other, denoted by $E
{}_2\!\!\!\sim F$, if there is a measurable partition $\{ E_n\}$ of
$E$, such that $\{ 2^n E_n \}$ is a measurable partition of $F$. A
measurable set $E$ is called a {\it $2\pi$-translation generator of
a measurable partition of $\mathbb{R}$} if $\{ E + 2n\pi \}_{n\in
\mathbb{Z}}$ forms a measurable partition of $\mathbb{R}$.
Similarly, a measurable set $F$ is called a {\it $2$-dilation
generator of a measurable partition of $\mathbb{R}$} if $\{ 2^n F
\}_{n\in \mathbb{Z}}$ forms a measurable partition of $\mathbb{R}$.

Lemma 4.3 in [DL] gives the following characterization of wavelet
sets, which was also obtained independently in [FW] using different
techniques.

{\it Let $E\subseteq \mathbb{R}$ be a measurable set. Then
$E$ is a wavelet set if and only if $E$ is both a $2\pi$-translation
generator of a measurable partition of $\mathbb{R}$ and a $2$-dilation
generator of a measurable partition of $\mathbb{R}$. }

Again from [DL], suppose that $E, F$ are wavelet sets, and let
$\sigma : E\rightarrow F$ be the bijective map (modulo null sets)
implementing the $2\pi$-translation congruence. Then $\sigma$ can be
extended to a bijective (modulo null sets) measurable map on
$\mathbb{R}$ by defining $\sigma(0)=0$ and $\sigma(s) = 2^n
\sigma(2^{-n}s)$ for each $s\in 2^n E, n\in \mathbb{Z}$. This map is
denoted by $\sigma_E^F$ and called the {\it interpolation map} for
the ordered pair of wavelet sets $(E, F)$.

An ordered pair of wavelet sets $(E, F)$ is called an {\it
interpolation pair} if $\sigma_E^F \circ \sigma_E^F :=
{(\sigma_E^F)}^2 = id_{\mathbb{R}}$. In general, an {\it
interpolation family of wavelet sets} is a family $\mathcal{F}$ of
wavelet sets such that $\{ \, \sigma_E^F \, |\, F\in \mathcal{F}
\,\}$ is a group under composition of maps for some $E\in
\mathcal{F}$.
\\
\\
{\bf Question A.} {\it Let $E, F$ be wavelet sets. Does $\sigma_E^F
(E \cup F) \subseteq E\cup F$ imply that $(\sigma_E^F)^2 =
id_{\mathbb{R}}$? }
\\

We answer Question A affirmatively. We then give an elementary
example that shows that this result need not hold for a {\it triple}
of wavelet sets; However, there is a natural modification which does
make sense for $n$-tuples, and we pose it as an open question
(Question C).

Given a pair of wavelet sets, it is not obvious at all upon initial
inspection whether they actually form an interpolation pair. And
constructing interpolation pairs can be hard. A basic problem
from [DL] that still remains open is the question: Given an
arbitrary dyadic wavelet set $E$ in $\mathbb{R}^1$, is there necessarily a
second distinct wavelet set $F$ such that $(E,F)$ is an
interpolation pair? Partly to address this problem, and partly for
intrinsic interest, Han introduced and studied properties of
congruence domains in [Han1] and [Han2]. Given a pair of wavelet sets
$(E, F)$, the {\it domain of $2\pi$-congruence} of $\sigma_E^F$,
denoted by $\mathcal{D}_E^F$, is the set of all points $s\in
\mathbb{R}$ such that $\sigma_E^F(s) - s$ is an integral multiple of
$2\pi$. There is a close relation between this and the interpolation
map of a pair of wavelet sets. Han asked the following question.
\\
\\
{\bf Question B.} {\it Let $E, F$ be wavelet sets. Does
$\mathcal{D}_E^F = \mathcal{D}_F^E$ imply that $(\sigma_E^F)^2 =
id_{\mathbb{R}}$? }
\\

If $(E,F)$ is an interpolation pair, then it is easily verified that
the domains of $2\pi$-congruence of $(E,F)$ and $(F,E)$ are the
same. So the above question just asks if the converse is true. In
many cases it is true.  But it is not universally true.  We answer
Question B negatively by constructing a counterexample. The key is a
special lemma on constructing wavelet sets, which is also used to
build examples for several other questions.

\section{Solution to question A}

The following theorem provides an answer to Question A.
\\
\\
{\bf Theorem 1.} {\it Let $(E, F)$ be wavelet sets. The following
two statements are equivalent:

(i) $(E, F)$ is an interpolation pair,

(ii) $\sigma_E^F (E\cup F) \subseteq E\cup F$.}
\\
\\
{\it Proof.} (i) $\Rightarrow$ (ii). Suppose that $(E, F)$ is an
interpolation pair. Then $(\sigma_E^F)^2 = id_{\mathbb{R}}$. Since
$(\sigma_E^F)^{-1} = \sigma_F^E$, it follows that
$\sigma_E^F = \sigma_F^E$. Observe that
$\sigma_E^F (E) = F$, $\sigma_F^E (F \backslash E) =
E\backslash F$, and so
\[ \sigma_E^F (E\cup F)  =  \sigma_E^F(E)
\cup \sigma_E^F(F \backslash E) =  F \cup \sigma_F^E(F\backslash
E) = F \cup E\backslash F = E \cup F. \]

(ii) $\Rightarrow$ (i). Suppose that $\sigma_E^F(E\cup F) \subseteq
E\cup F$. Since $\sigma_E^F(E) = F$ and $\sigma_E^F$ is bijective
(modulo null sets), we must have $\sigma_E^F(F\backslash E)
\subseteq E\backslash F$. We will prove this by way of
contradiction, and a diagram is included at the end of the proof to
help in navigating between $E$ and $F$, with arrows representing
$\sigma^F_E$.

Assume that $(E, F)$ is not an
interpolation
pair. Let $G_0 = \{\, s\in E \, |\, (\sigma_E^F)^2 (s) \ne s\,\}$. Then
$G_0 \subseteq E\backslash F$ is Lebesgue measurable and has positive
measure. Since $\{\, G_0 \cap (F-2n\pi) \,\}_{n\in \mathbb{Z}}$ forms a
measurable partition of $G_0$, it follows that $G_0 \cap (F-2n_1\pi)$ has
positive measure for some $n_1 \in
\mathbb{Z}$. Denote this set by $G_1$.  Then $\sigma_E^F (G_1) =
G_1 + 2n_1\pi \subseteq F\backslash E$.  Since $E\backslash F$ is
$2$-dilation congruent to $F\backslash E$, following the similar
discussion, there exists a measurable subset $G_2$ of $G_1$ with positive
measure, such that $2^{-k_1}(G_2+2n_1\pi) \subseteq E\backslash F$
for some integer $k_1$. Then, there exists a measurable subset $G_3$ of
$G_2$ with positive measure such that $\sigma_E^F (2^{-k_1}(G_3+2n_1\pi)) =
2^{-k_1}(G_3+2n_1\pi) + 2n_2\pi$ for some integer $n_2$. Thus,
\begin{eqnarray*}
(\sigma_E^F)^2 (G_3) & = & \sigma_E^F (G_3 +2n_1\pi) = 2^{k_1} \cdot
\sigma_E^F ( 2^{-k_1}(G_3 +2n_1\pi) ) \\
& = & 2^{k_1} \cdot ( 2^{-k_1}(G_3 +2n_1\pi) + 2n_2\pi) =
G_3 +2n_1\pi + 2^{k_1} \cdot 2n_2\pi \subseteq E\backslash F .
\end{eqnarray*}
Since both $G_3, G_3 +2n_1\pi +2^{k_1} \cdot 2n_2\pi \subseteq E\backslash
F$ and they are distinct by assumption, we must have $n_1 + 2^{k_1}\cdot
n_2 \notin \mathbb{Z}$, which implies that $k_1 \le
-1$. Similarly, there exists a measurable subset $G_4$
of $G_3$ with positive measure, such that $\sigma_E^F (G_4+2n_1\pi
+2^{k_1}\cdot 2n_2\pi) = G_4 + 2n_1\pi + 2^{k_1}\cdot 2n_2\pi + 2n_3\pi$
for some integer $n_3$. Then,
\begin{eqnarray*}
\lefteqn{(\sigma_E^F)^2 (2^{-k_1} (G_4 +2n_1\pi)) =
 \sigma_E^F (2^{-k_1} (G_4 +2n_1\pi) + 2n_2\pi) } \\
&=& 2^{-k_1} \cdot \sigma_E^F (G_4 + 2n_1\pi +
2^{k_1} \cdot 2n_2\pi)
= 2^{-k_1} \cdot (G_4 + 2n_1\pi + 2^{k_1} \cdot 2n_2\pi + 2n_3\pi )\\
& =&
2^{-k_1} (G_4 +2n_1\pi) + 2n_2\pi + 2^{-k_1}\cdot 2n_3\pi \subseteq
E\backslash F .
\end{eqnarray*}
Since both $2^{-k_1} (G_4 +2n_1\pi)$ and $
2^{-k_1} (G_4 +2n_1\pi) + 2n_2\pi + 2^{-k_1}\cdot 2n_3\pi$ are contained
in $E\backslash F$ and $n_2 + 2^{-k_1}\cdot n_3$
must be an integer, we must have $n_2  + 2^{-k_1}\cdot n_3 = 0$. Then,
since $\sigma_E^F$ maps both $G_4$ and $G_4 + 2n_1\pi +2^{k_1}\cdot
2n_2\pi$ to $G_4 + 2n_1\pi$, $G_4$ and $G_4+2n_1\pi+2^{k_1}\cdot 2n_2\pi$
must be the same, then $(\sigma_E^F)^2 (s) = s,
\forall s \in G_4 \subseteq G_0$, contradicting the assumption. Thus,
$(\sigma_E^F)^2 = id_{\mathbb{R}}$ almost everywhere.
\begin{center}
\begin{picture}(400, 95)(0,0)
\put(0,90){$E\backslash F$ :}
\put(0,60){$F\backslash E$ :}
\put(0,30){$E\backslash F$ :}
\put(0,0){$F\backslash E$ :}

\put(105,90){$G_0$}
\put(110,83){\vector(0,-1){10}}
\put(90,60){$G_1+2n_1\pi$}
\put(110,53){\vector(0,-1){10}}
\put(60,30){$G_3+2n_1\pi +2^{k_1}\cdot 2n_2\pi$}
\put(110,23){\vector(0,-1){10}}
\put(40,0){$G_4+2n_1\pi +2^{k_1}\cdot 2n_2\pi + 2n_3\pi$}

\put(255,90){$2^{-k_1}(G_2+2n_1\pi)$}
\put(288,83){\vector(0,-1){10}}
\put(230,60){$2^{-k_1}(G_3+2n_1\pi)+2n_2\pi$}
\put(288,53){\vector(0,-1){10}}
\put(200,30){$2^{-k_1}(G_4+2n_1\pi)+2n_2\pi+2^{-k_1}\cdot 2n_3\pi$}
\put(288,23){\vector(0,-1){10}}
\put(283,0){$\cdots$}
\end{picture}
\end{center}
$\ \ \ \ \ \ \ \ \ \ \ \ \ \ \ \ \ \ \ \ \ \ \ \ \ \ \ \ \ \ \ \ \
\ \ \ \ \ \ \ \ \ \ \ \ \ \ \ \ \ \ \ \ \ \ \ \ \ \ \ \ \ \ \
\ \ \ \ \ \ \ \ \ \ \ \ \ \ \ \ \ \ \ \ \ \ \ \ \ \ \ \ \ \ \
\ \ \ \ \ \ \ \ \ \ \ \ \ \ \ \ \ \ \ \ \ \ \ \ \ \ \ \ \ \ \ \square$
\\

Let $\mathcal{F}$ be a general interpolation family of wavelet sets,
and fix $E_1\in \mathcal{F}$.  Then for arbitrary $E_2, E_3 \in
\mathcal{F}$, we have $\sigma_{E_1}^{E_2} \circ \sigma_{E_1}^{E_3} =
\sigma_{E_1}^{E_4}$ for some $E_4 \in \mathcal{F}$.  Thus  $
\sigma_{E_1}^{E_2}(E_3) = \sigma_{E_1}^{E_2} \circ
\sigma_{E_1}^{E_3} (E_1) = \sigma_{E_1}^{E_4}(E_1) = E_4$.  It
follows that  $\sigma_{E_1}^{E_{\lambda}} (\bigcup_{E\in
\mathcal{F}} E) \subseteq \bigcup_{E\in \mathcal{F}} E$ for each
$E_{\lambda} \in \mathcal{F}$. However, the converse may be false as
the following example shows.
\\
\\
{\bf Example 2.} Let $(E, F)$ be an interpolation pair of wavelet
sets and suppose  $G \subset E \cup F$ is a wavelet set contained in
the union which is distinct from $E$ and $F$. Observe that
$E\backslash F$ is both $2\pi$-translation congruent and
$2$-dilation congruent to $F\backslash E$, thus $G \subset E \cup F$
must contain $E\cap F$ to be a wavelet set. Since $\sigma_E^G =
id_{\mathbb{R}}$ on $E\cap G$ and $\sigma_E^G = \sigma_E^F$ on
$E\backslash G$, $(\sigma_E^G)^2 = id_{\mathbb{R}}$. Similarly,
$(\sigma_F^G)^2 = id_{\mathbb{R}}$. Then, $\sigma_E^F \circ
\sigma_E^G = \sigma_F^G \ne id_{\mathbb{R}}$, which implies that $\{
id_{\mathbb{R}}, \sigma_E^F ,\sigma_E^G \}$ is not a group under
composition of maps. Hence $\{ E,F,G\}$ is not an interpolation
family. However, notice that $\{ E\backslash G , G\backslash F ,
E\cap F , G\backslash E, F\backslash G \}$ forms a measurable
partition of $E\cup F$. $E\backslash G$ is both $2\pi$-translation
and $2$-dilation congruent to $G\backslash E$ and $F\backslash G$ is
both $2\pi$-translation and $2$-dilation congruent to $G\backslash
F$ which is contained in $E\cap G$. We have $\sigma_E^F (E\cup F
\cup G) = \sigma_E^F(E\cup F) = E\cup F$, and
\begin{eqnarray*}
\sigma_E^G (E\cup F \cup G) &=& \sigma_E^G (E\cup F) =
\sigma_E^G(E) \cup \sigma_E^G(F\backslash E) \\
& = & G \cup \sigma_E^G
(G\backslash E) \cup \sigma_E^G (F\backslash G)  =  G \cup  E\backslash G
\cup F\backslash G = E\cup F.
\end{eqnarray*}
\ \ \ \ \ \ \ \ \ \ \ \ \ \ \ \ \ \ \ \ \ \ \ \ \ \ \ \ \ \ \ \ \
\ \ \ \ \ \ \ \ \ \ \ \ \ \ \ \ \ \ \ \ \ \ \ \ \ \ \ \ \ \ \ \ \
\ \ \ \ \ \ \ \ \ \ \ \ \ \ \ \ \ \ \ \ \ \ \ \ \ \ \ \ \ \ \ \ \
\ \ \ \ \ \ \ \ \ \ \ \ \ \ \
$\square$
\\

Observe that in the above example, if let $H = E\backslash G \cup
(E\cap F) \cup F\backslash G$, then $H\subseteq E\cup F$ is also a
wavelet set and $\sigma_E^H = \sigma_E^F \circ \sigma_E^G =
\sigma_E^G \circ \sigma_E^F$. Then $\{ E,F,G, H\}$ is an
interpolation family since $\{ \sigma_E^E=id_{\mathbb{R}},
\sigma_E^F, \sigma_E^G, \sigma_E^H \}$ is actually isomorphic to the
Klein four group. This motivates the following question which seems
difficult and which we pose as an open problem:
\\
\\
{\bf Question C.} Let $\mathcal{F}$ be a finite collection of
wavelet sets. Fix $E_1\in \mathcal{F}$. Suppose
$\sigma_{E_1}^{E_{\lambda}} (\bigcup_{E\in \mathcal{F}} E)$ is
contained in $\bigcup_{E\in \mathcal{F}} E$ for each $E_{\lambda}\in
\mathcal{F}$.  Is  $\mathcal{F}$ a {\it{subset}} of a finite
interpolation family of wavelet sets?  That is, can $\mathcal{F}$
always be {\it{extended}} to a finite interpolation family?

\section{A counterexample to question B}

Before giving a counterexample to Question B, we present a
characterization of interpolation pairs of wavelet sets and equality of
$2\pi$-congruence domains of associated interpolation maps. This is
actually the motivation that led us to find a counterexample.

Given a pair of wavelet sets $(E, F)$, there exists a measurable partition
$\{ E_n\}_{n\in \mathbb{Z}}$ of $E$ such that $\{ E_n +2n\pi \}_{n\in
\mathbb{Z}}$ is a measurable partition of $F$. Similarly, there exists
another measurable partition $\{ E^k \}_{k\in \mathbb{Z}}$ of $E$ such
that $\{ 2^k E^k \}_{k\in \mathbb{Z}}$ is another measurable partition of
$F$. Denote $E_n \cap E^k$ by $E_{n,k}$, then $\{ E_{n,k} \}_{n,k\in
\mathbb{Z}}$ is a measurable partition of $E$, and both $\{ E_{n,k} +
2n\pi \}_{n,k\in \mathbb{Z}}$ and $\{ 2^k E_{n,k} \}_{n,k\in \mathbb{Z}}$
are measurable partitions of $F$. Observe that if $E_{n,k}$ has positive
measure, then $n=0$ whenever $k=0$, and $E_{n,k} \cap (2^l E_{m,l} -
2n\pi), n,k,m,l \in \mathbb{Z}\backslash 0$ forms a partition of
$E\backslash F$. \\
\\
{\bf Theorem 3.} {\it Let $E,F$ be wavelet sets, and let $\{
E_{n,k}\}_{n,k\in\mathbb{Z}}$ be as defined above. Then:

(i) $(E, F)$ is an interpolation pair if and only
if $E_{n,k} \cap (2^l E_{m,l} - 2n\pi)$ with positive measure for some
nonzero integers $n,k,m,l$ implies that $n + 2^l\cdot m
=0$. Furthermore, for such a set, we also have $k=-l$.

(ii) $\mathcal{D}_E^F = \mathcal{D}_F^E$ if and only if  $E_{n,k} \cap
(2^l E_{m,l} - 2n\pi)$ with positive measure for some nonzero integers
$n,k,m,l$ implies that $[n] + l = [m]$. Here, $[n]$ denotes
the smallest integer $k$ such that $2^k\cdot n\in \mathbb{Z}$. }
\\
\\
{\it Proof.} (i) Suppose that $(E, F)$ is an interpolation pair.
Then $(\sigma_E^F)^2 = id$. Since
\begin{eqnarray*}
(\sigma_E^F)^2 (E_{n,k} \cap (2^lE_{m,l}-2n\pi))
& = & \sigma_E^F ( (E_{n,k} + 2n\pi) \cap 2^l E_{m,l} ) \\
& = & 2^l \cdot \sigma_E^F ( 2^{-l}(E_{n,k} + 2n\pi) \cap E_{m,l} ) \\
& = & 2^l\cdot (\, ( 2^{-l}(E_{n,k} + 2n\pi) \cap E_{m,l} ) +2m\pi ) \\
& = &  (E_{n,k} + 2n\pi + 2^l\cdot 2m\pi) \cap 2^l (E_{m,l}+2m\pi),
\end{eqnarray*}
it follows that if $E_{n,k} \cap (2^l E_{m,l} -2n\pi)$ has positive
measure, then $2n\pi +2^l\cdot 2m\pi =0$, i.e., $n + 2^l\cdot m =0$.

Conversely, suppose that $E_{n,k} \cap (2^l E_{m,l} -2n\pi),
n,k,m,l\ne 0$ having positive measure implies that $n + 2^l\cdot m
=0$. Observe that $\{ E_{n,k} \cap (2^l E_{m,l} -2n\pi)
\}_{n,k,m,l\ne 0}$ forms a measurable partition of $E\backslash F$.
Then, a similar argument shows that $(\sigma_E^F)^2 =
id_{\mathbb{R}}$.

Furthermore, if $E_{n,k} \cap (2^l E_{m,l} -2n\pi)$ having positive
measure implies that $n+2^l\cdot m =0$, then $( E_{m,l} + 2m\pi)
\cap 2^{-l} E_{n,k} = 2^{-l}\cdot ( E_{n,k} \cap (2^l E_{m,l}
-2n\pi)\,)$ also has positive measure. Since $E_{m,l} +2m\pi
\subseteq F$ and $2^{-l} E_{n,k} \subseteq F$ only if $k=-l$, we
must have $k=-l$.

(ii) Observe that
\begin{eqnarray*}
\mathcal{D}_E^F &=& \dot{\bigcup_{n,k\ne 0}} (
\dot{\bigcup_{j\ge [n]}} 2^j E_{n,k} )  \,\,\, \dot{\cup} \,\,\,
 \dot{\bigcup_{j}} 2^j (E\cap F)  \\
&=& \dot{\bigcup_{n,k,m\ne 0}} ( \dot{\bigcup_{j\ge [n]}}
2^j (E_{n,k}\cap 2^{-k}(E_m + 2m\pi)) ) \,\,\,
\dot{\cup} \,\,\,  \dot{\bigcup_{j}} 2^j (E\cap F)  ,
\end{eqnarray*}
and
\begin{eqnarray*}
\mathcal{D}_F^E &=&  \dot{\bigcup_{n,k\ne 0}} (
\dot{\bigcup_{j\ge [n]}} 2^j (E_{n,k}+2n\pi) ) \,\,\,\dot{\cup}
\,\,\, \dot{\bigcup_{j}} 2^j (E\cap F)  \\
&=&  \dot{\bigcup_{n,k,m,l\ne  0}} ( \dot{\bigcup_{j\ge [n]}}
2^j ( (E_{n,k}+2n\pi) \cap 2^lE_{m,l}) ) \,\,\,
\dot{\cup} \,\,\,  \dot{\bigcup_{j}} 2^j (E\cap F)  \\
&=& \dot{\bigcup_{n,k,m,l\ne 0}} ( \dot{\bigcup_{j\ge [n]}}
2^{j+l} ( E_{m,l} \cap 2^{-l}
(E_{n,k}+2n\pi) ) ) \,\,\,\dot{\cup} \,\,\,
\dot{\bigcup_{j}} 2^j (E\cap F)   \\
&=& \dot{\bigcup_{n,k,m,l\ne 0}} ( \dot{\bigcup_{j\ge [m]}}
2^{j+k} ( E_{n,k} \cap 2^{-k} (E_{m,l}+2m\pi) ) ) \,\,\,
\dot{\cup} \,\,\,  \dot{\bigcup_{j}} 2^j (E\cap F) \\
&=& \dot{\bigcup_{n,k,m\ne 0}} ( \dot{\bigcup_{j\ge [m]}} 2^{j+k}
(E_{n,k} \cap 2^{-k}( E_m + 2m\pi))) \,\,\, \dot{\cup} \,\,\,
 \dot{\bigcup_{j}} 2^j (E\cap F) ,
\end{eqnarray*}
where $\dot{\cup}$ denotes disjoint union. The next to last equality
comes from changing indices. Thus, $\mathcal{D}_E^F =
\mathcal{D}_F^E$ implies that $[m] + k = [n]$, if $E_{n,k} \cap
2^{-k} (E_m + 2m\pi)$ has positive measure, or equivalently, if
$E_{m,l} \cap (2^k E_{n,k} - 2m\pi)$ has positive measure for some
$l\neq 0$, since $2^k ( E_{n,k} \cap 2^{-k} (E_m + 2m\pi) ) - 2m \pi
= (2^k E_{n,k} -2m\pi ) \cap E_m = (2^k E_{n,k} -2m\pi ) \cap (
\dot{\bigcup}_lE_{m,l} )$ . The converse direction can be shown by
reversing the above discussion. \ \ \ \ \ \ \ \ \ \ \ \ \ \ \ \ \ \
\ \ \ \ \ \ \ \ \ \ \ \ \ \ \ \ \ \ \ \ \ \ \ \ \ \ \ \ \ \ \ \ \ \
\ \ \ \ \ \ \ \ \ \ \ \ \ \ \ \ \ \ \ \ \ \ \ \ \ \ \ \ \ \ \ \ \ \
\ \ \ \ \ \ \ $\square$
\\

The basic idea we use is to construct two certain measurable subsets
of $\mathbb{R}$, which are both $2\pi$-translation and $2$-dilation
congruent to each other, yet the interpolation map restricted
to the union of $2$-dilates of which does not have the property that
its square equals the identity map, and then to construct the
remaining pieces of wavelet sets for these two sets so that the
congruency domains match up. This type of approach was used by D.
Speegle in [S1] and Q. Gu in [Gu1], and a necessary and sufficient
condition has also been given for a measurable set being contained
in some wavelet set in [IP]. However, the methods and the
constructions we use in the present paper are completely different
from those in [S1], [Gu1] and[IP].
\\
\\
{\bf Theorem 4.} {\it The answer to Question B is no.}
\\

Before proving Theorem 4, we require a technical lemma. This will also be
used in constructing several other examples in the next section.
\\
\\
{\bf Lemma 5.}  {\it Let $E, F\subseteq {\mathbb{R}}^n$ be
Lebesgue measurable sets with finite positive measure. If there exist
$n_1, n_2, k_1, k_2 \in {\mathbb{Z}}^{(n)}$, such that modulo null sets,
\[ 2^{k_1}F \subseteq E + 2n_1\pi \mbox{\, \, \, and \, \, \, }
   E + 2n_2\pi \subseteq 2^{k_2}F , \]
and $(E+2n_2\pi) \cap 2^{k_1} F$ is a null set, then there exists a
Lebesgue measurable set $G \subseteq (E+2n_2\pi) \cup 2^{k_1} F$
such that $G$ is both $2\pi$-translation congruent to $E$ and
$2$-dilation congruent to $F$. In fact $G$ can be taken as:
\begin{eqnarray*}
G &=& \bigcup_{i=0}^{\infty} S^i ( 2^{k_1}F \backslash 2^{k_1-k_2} (E +
2n_2\pi) ) \\
  &  & \cup \,\,\, (E+2n_2\pi) \backslash ( \bigcup_{i=1}^{\infty}
2^{k_2-k_1} \cdot S^i
(2^{k_1}F \backslash 2^{k_1-k_2} (E + 2n_2\pi) )) ,
\end{eqnarray*}
where $S(x) = 2^{k_1-k_2} (x+2(n_2-n_1)\pi), \forall x\in
{\mathbb{R}}^n$. }
\\
\\
{\it Proof.} By hypothesis,
\[ 2^{k_1-k_2} ( E+2n_2\pi) \subseteq 2^{k_1} F \subseteq E+2n_1\pi. \]
Since $\mu(2^{k_1-k_2} (E+2n_2\pi)) = 2^{k_1-k_2}\mu(E) \le
\mu(E+2n_1\pi) = \mu(E)$, $k_1-k_2 \le0$. If $k_1 = k_2$, then $E+2n_2\pi
\subseteq 2^{k_1}F$, contradicts the hypothesis that
$(E+2n_2\pi) \cap
2^{k_1}F$ is a null set. Thus, $k_1 - k_2 <0$. Construct a sequence of
measurable sets $\{
G_i\}_{i\in \mathbb{N}}$ as follows. Let
\begin{eqnarray*}
G_0 & = &  2^{k_1}F \,\, \backslash \,\,  2^{k_1-k_2} (E+2n_2\pi) \\
 & = & 2^{k_1}F \,\, \backslash \,\, 2^{k_1-k_2} (E+2n_1\pi +
2(n_2-n_1)\pi) \subseteq  (E+2n_1\pi) \,\, \backslash \,\,
S(E+2n_1\pi) .
\end{eqnarray*}
Let $G_i = S^i(G_0) \subseteq S^i(E+2n_1\pi) \backslash
S^{i+1}(E+2n_1\pi)$, for each $i\in \mathbb{N}$. Notice that the
measure of $G_i$ is bounded by $2^{i\cdot (k_1-k_2)} \cdot \mu(E)$,
which will approach $0$ as $i$ approaches infinity. Furthermore,
$S(E+2n_1\pi) \subseteq E+2n_1\pi$, and it follows that the $G_i$'s
are disjoint and $\bigcup_{i=0}^{\infty} G_i \subseteq 2^{k_1} F$.
Also by definition, $\cup_{i=1}^{\infty} G_i \subseteq S(E+2n_1\pi)$
and $2^{k_2-k_1} G_{i+1} = G_i + 2(n_2-n_1)\pi$, $\forall i\ge 0$.
Let
\[ G = ( \bigcup_{i=0}^{\infty} G_i ) \,\,\, \dot{\cup} \,\,\,
(E+2n_2\pi) \backslash (\bigcup_{i=1}^{\infty} 2^{k_2-k_1} G_i) . \]
Then, since $2^{k_1 - k_2} ( (E+2n_2\pi) \backslash
(\bigcup_{i=1}^{\infty} 2^{k_2-k_1} G_i) ) = S(E+2n_1\pi) \backslash
(\bigcup_{i=1}^{\infty} G_i)$, and
$\{ G_i\} _{i=0}^{\infty} \cup \{ S(E+2n_1\pi)
\backslash ( \bigcup _{i=1}^{\infty} G_i ) \}$ constitutes a measurable
partition of $2^{k_1}F$, it's clear that $G$ is $2$-dilation congruent
to $F$. On the other hand, since
\[ \bigcup _{i=1}^{\infty} 2^{k_2-k_1} G_i
     = \bigcup _{i=0}^{\infty} (G_i + 2(n_2-n_1)\pi) , \]
$G$ is $2\pi$-translation congruent to $E$.
\ \ \ \ \ \ \ \ \ \ \ \ \ \ \ \ \ \ \ \ \ \ \ \ \ \ \ \ \ \ \ \ \ \ \ \
\ \ \ \ \ \ \ \ \ \ \ \ \ \ \ \ \ \ \ \ \ \ \ \ \ \ \ \ \ \ \ \ \ \ \ \
\ \ \ \ \ \ \ \ \ \ \ \
$\square$
\\
\\
{\it Proof of Theorem 4.} Let
$E_1 = [\frac{33}{16}\pi, \frac{34}{16}\pi),
E_2 = [\frac{33}{8}\pi, \frac{34}{8}\pi) + 12\pi,
E_3 = [\frac{33}{4}\pi, \frac{34}{4}\pi) + 8\pi$ and $
E_4 = [\frac{33}{2}\pi, \frac{34}{2}\pi) + 96\pi$. Let
$F_1 = [\frac{33}{4}\pi, \frac{34}{4}\pi),
F_2 = [\frac{33}{16}\pi, \frac{34}{16}\pi) + 6\pi,
F_3 = [\frac{33}{2}\pi, \frac{34}{2}\pi) + 16\pi$ and $
F_4 = [\frac{33}{8}\pi, \frac{34}{8}\pi) + 24\pi$.
Notice that $E_1, E_2, E_3, E_4$ are both $2\pi$-translation
congruent and $2$-dilation congruent to disjoint pieces of $[-2\pi,
-\pi) \cup [\pi, 2\pi)$, and so are $F_1, F_2, F_3,
F_4$. Furthermore, \begin{eqnarray*}
E_1\cup E_2\cup E_3\cup E_4 & \sim_{2\pi} &  F_1\cup
F_2\cup F_3\cup F_4 \ \sim_{2\pi} \ [\frac{\pi}{16}, \pi), \\
E_1\cup E_2\cup E_3\cup E_4  & {}_2\!\!\!\sim &  F_1\cup
F_2\cup F_3\cup F_4 \ {}_2\!\!\!\sim \\
&   & [\pi +\frac{\pi}{128}, \pi + \frac{\pi}{16}) \cup
[\frac{7}{4}\pi +\frac{\pi}{128}, \frac{7}{4}\pi +\frac{\pi}{64}).
\end{eqnarray*}
Let $G$ be the measurable set determined by Lemma 5
and the following two containment relations:
\[ \frac{1}{32} (  [\pi, \pi +\frac{\pi}{128}) \cup [\pi +
\frac{\pi}{16}, \frac{7}{4}\pi +\frac{\pi}{128}) \cup [\frac{7}{4}\pi
+\frac{\pi}{64}, 2\pi) ) \ \subseteq \ [0, \frac{\pi}{16}), \]
\[ [0, \frac{\pi}{16}) + 8\pi \ \subseteq \ 8\cdot ( [\pi, \pi
+\frac{\pi}{128}) \cup [\pi +
\frac{\pi}{16}, \frac{7}{4}\pi +\frac{\pi}{128}) \cup [\frac{7}{4}\pi
+\frac{\pi}{64}, 2\pi) ). \]
Then, $G$ is both $2\pi$-translation congruent to $[0, \frac{\pi}{16})$
and $2$-dilation congruent to  $[\pi, \pi +\frac{\pi}{128}) \cup [\pi +
\frac{\pi}{16}, \frac{7}{4}\pi +\frac{\pi}{128}) \cup [\frac{7}{4}\pi
+\frac{\pi}{64}, 2\pi)$. Hence,
\[ E := [-2\pi, -\pi) \cup G \cup E_1 \cup E_2 \cup E_3\cup E_4, \]
\[ F := [-2\pi, -\pi) \cup G \cup F_1 \cup F_2 \cup F_3\cup F_4 \]
are both wavelet sets. The associated interpolation map $\sigma_E^F$ is
defined by:
\[ \sigma_E^F(s) = \left\{ \begin{array}{ll}
       s+6\pi & \mbox{ if } s\in E_1 \\
       s+12\pi & \mbox{ if } s\in E_2 \\
       s-8\pi & \mbox{ if } s\in E_3 \\
       s-80\pi & \mbox{ if } s\in E_4 \\
       s  & \mbox{ if } s\in [-2\pi, -\pi) \cup G .
                           \end{array}
                   \right. \]
Straightforward computation shows that
\[ \mathcal{D}_E^F = = (-\infty, 0) \cup
\bigcup_{k\in \mathbb{Z}} 2^k G  \cup  \bigcup_{k\ge 0} (E_1 \cup
\frac{1}{2} E_2 \cup \frac{1}{4} E_3 \cup \frac{1}{8} E_4 ) \]
and
\[ \mathcal{D}_F^E = = (-\infty, 0) \cup
\bigcup_{k\in \mathbb{Z}} 2^k G  \cup  \bigcup_{k\ge 0} (F_2 \cup
\frac{1}{2} F_4 \cup \frac{1}{4} F_1 \cup \frac{1}{8} F_3 ) . \]
Observe that $F_2 = \frac{1}{2} E_2,
\frac{1}{2} F_4 = \frac{1}{8} E_4,
\frac{1}{4} F_1 = E_1$, and $
\frac{1}{8} F_3 = \frac{1}{4} E_3$, which implies that $\mathcal{D}_E^F =
\mathcal{D}_F^E$.
However, since
\[ (\sigma_E^F)^2 (E_1) = \sigma_E^F (E_1 + 6\pi) = \sigma_E^F
(\frac{1}{2}E_2) = \frac{1}{2} \sigma_E^F(E_2) = \frac{1}{2} (E_2 +
12\pi) = E_1 + 12\pi, \]
$(E,F)$ is not an interpolation pair.
\ \ \ \ \ \ \ \ \ \ \ \ \ \ \ \ \ \ \ \ \ \ \ \ \ \ \ \ \ \ \ \ \ \ \ \ \
\ \ \ \ \ \ \ \ \ \ \ \ \ \ \ \ \ \ \ \ \ \ \ \ \ \ \ \ \ \ \ \ \ \ \ \ \
\ \ \ \ \ \ \ \ \ \
$\square$
\\

It turns out that if one of $E, F$ is the Shannon set, then the equality
of integral domains does imply that $(E, F)$ is an interpolation pair.
\\
\\
{\bf Proposition 6.} {\it Let $E,F$ be wavelet sets. If $E = [-2\pi,
-\pi) \cup [\pi, 2\pi)$, then $\mathcal{D}_E^F = \mathcal{D}_F^E$
implies that $(E,F)$ is an interpolation pair. }
\\
\\
{\it Proof.} Let $E_n = E\cap (F-2n\pi)$, $ n\in \mathbb{Z}$. Then
$\{ E_n \}_{n\in \mathbb{Z}}$ is a measurable partition of $E$ and
$\{ E_n+2n\pi\}$ is a measurable partition of $F$. Let $E_n^- = E_n
\cap [-2\pi, -\pi), E_n^+ = E_n \cap [\pi, 2\pi)$. Then, we have the
following diagram.

\begin{center}
\begin{picture}(400,80)(0,0)
\put(138,75){$E_n^-$}
\put(258,75){$E_n^+$}

\put(0,70){$E:$}
\put(143,70){\oval(40,5)[b]}
\put(183,70){\oval(40,5)[b]}
\put(223,70){\oval(40,5)[b]}
\put(263,70){\oval(40,5)[b]}

\put(118,55){$-2\pi$}
\put(158,55){$-\pi$}
\put(203,55){$0$}
\put(243,55){$\pi$}
\put(278,55){$2\pi$}

\put(26,30){$\cdots$}
\put(45,30){$E_{-1}^-\!\!-\!\!2\pi$}
\put(85,30){$E_{-2}^+\!\!-\!\!4\pi$}
\put(138,30){$E_0^-$}
\put(165,30){$E_{-1}^+\!\!-\!\!2\pi$}
\put(205,30){$E_1^-\!\!+\!\!2\pi$}
\put(258,30){$E_0^+$}
\put(285,30){$E_2^-\!\!+\!\!4\pi$}
\put(325,30){$E_1^+\!\!+\!\!2\pi$}
\put(366,30){$\cdots$}

\put(0,25){$F:$}
\put(23,25){\oval(40,5)[br]}
\put(63,25){\oval(40,5)[b]}
\put(103,25){\oval(40,5)[b]}
\put(143,25){\oval(40,5)[b]}
\put(183,25){\oval(40,5)[b]}
\put(223,25){\oval(40,5)[b]}
\put(263,25){\oval(40,5)[b]}
\put(303,25){\oval(40,5)[b]}
\put(343,25){\oval(40,5)[b]}
\put(383,25){\oval(40,5)[bl]}

\put(38,10){$-4\pi$}
\put(78,10){$-3\pi$}
\put(118,10){$-2\pi$}
\put(158,10){$-\pi$}
\put(203,10){$0$}
\put(243,10){$\pi$}
\put(278,10){$2\pi$}
\put(318,10){$3\pi$}
\put(358,10){$4\pi$}

\end{picture}
\end{center}

Observe that
\[ \mathcal{D}_E^F = \dot{\bigcup_{n\ne 0, j\ge
[n]}} 2^j (E_n^- \cup E_n^+) \,\,\, \dot{\cup}\,\,\,  \dot{\bigcup_{j\in
\mathbb{Z}}} 2^j E_0 .
\]
Then
\begin{eqnarray*}
\mathcal{D}_E^F \cap [0,\pi) & =&  \dot{\bigcup_{n\ne 0
, [n] \le j < 0}} 2^j E_n^+ \,\,\,\dot{\cup}\,\,\,
 \dot{\bigcup_{j<0}} 2^j E_0^+  \\
&= &   \dot{\bigcup_{j<0}} 2^j E_0^+  \,\dot{\cup}\, \frac{1}{2} E_{\pm
2}^+ \,\dot{\cup}\, \frac{1}{2} E_{\pm 4}^+ \,\dot{\cup}\, \frac{1}{4}
E_{\pm 4}^+ \,\dot{\cup}\, \frac{1}{2} E_{\pm 6}^+ \,\dot{\cup}\,
\frac{1}{2} E_{\pm 8}^+ \,\dot{\cup}\, \frac{1}{4} E_{\pm 8}^+
\,\dot{\cup}\, \cdots,
\end{eqnarray*}
here $E_{\pm n} := E_n \cup E_{-n}$. $\mathcal{D}_E^F \cap
[\pi,2\pi) = [\pi, 2\pi)$. Similarly,
\[ \mathcal{D}_F^E =  \dot{\bigcup_{n \ne 0, j\ge
[n]}} 2^j ( (E_n^- \cup E_n^+)+2n\pi) \,\,\,\dot{\cup}\,\,\,
\dot{\bigcup_{j\in \mathbb{Z}}} 2^j E_0 . \]
Then, following from the diagram,
\[ \mathcal{D}_F^E \cap [0,\pi) = (E_1^- +2\pi) \,\,\, \dot{\cup}\,\,\,
  \dot{\bigcup_{j<0}} 2^j E_0^+  ,\]
\[ \mathcal{D}_F^E \cap [\pi,2\pi) = E_0^+ \,\dot{\cup}\, (
\dot{\bigcup_{j\ge 1}} 2^j(E_1^- +2\pi) \cap [\pi, 2\pi) ) \,\dot{\cup}\,
\frac{1}{2}(E_2^-+4\pi) \,\dot{\cup}\, \frac{1}{4} (E_4^- +
8\pi) \,\dot{\cup}\, \cdots. \]
Comparing $\mathcal{D}_E^F \cap [0,\pi)$ and $\mathcal{D}_F^E \cap
[0,\pi)$, we have
\begin{eqnarray}
E_1^- + 2\pi & = & \frac{1}{2} E_{\pm 2}^+ \,\dot{\cup}\, \frac{1}{2}
E_{\pm 4}^+ \,\dot{\cup}\,  \frac{1}{4} E_{\pm 4}^+ \,\dot{\cup}\,
\frac{1}{2} E_{\pm 6}^+ \,\dot{\cup}\,  \frac{1}{2} E_{\pm 8}^+
\,\dot{\cup}\, \frac{1}{4} E_{\pm 8}^+ \,\dot{\cup}\, \cdots.
\end{eqnarray}
Then, it is easy to see that $(\bigcup_{j\ge 1} 2^j (E_1^- +2\pi) ) \cap
[\pi, 2\pi) = \bigcup_{n\in 2\mathbb{Z}\backslash 0} E_n^+$.
Comparing $\mathcal{D}_E^F \cap [\pi, 2\pi)$ and $\mathcal{D}_F^E
\cap [\pi,2\pi)$, we have
\begin{eqnarray}
\dot{\bigcup_{n \in 2\mathbb{Z}+1}} E_n^+ & \subseteq & \frac{1}{2}
(E_2^- +4\pi) \,\dot{\cup}\, \frac{1}{4} (E_4^- + 8\pi) \,\dot{\cup}\,
\cdots .
\end{eqnarray}
By symmetry, we also have
\begin{eqnarray}
E_{-1}^+ - 2\pi & = & \frac{1}{2} E_{\pm 2}^- \,\dot{\cup}\, \frac{1}{2}
E_{\pm 4}^- \,\dot{\cup}\,  \frac{1}{4} E_{\pm 4}^- \,\dot{\cup}\,
\frac{1}{2} E_{\pm 6}^- \,\dot{\cup}\, \frac{1}{2} E_{\pm 8}^-
\,\dot{\cup}\, \frac{1}{4} E_{\pm 8}^- \,\dot{\cup}\, \cdots, \\
\dot{\bigcup_{n \in 2\mathbb{Z}+1}} E_n^- & \subseteq & \frac{1}{2}
(E_{-2}^+ -4\pi) \,\dot{\cup}\, \frac{1}{4} (E_{-4}^+ -
8\pi) \,\dot{\cup}\, \cdots .
\end{eqnarray}
(1) and (4) imply that $E_1^- +2\pi = \frac{1}{2} E_{-2}^+$, and
$E_n^-$ has measure zero for each odd integer $n$ except $1$, and
$E_n^+$ has measure zero for each even integer $n$ except $-2,0$.
(2) and (3) imply that  $E_{-1}^+ -2\pi = \frac{1}{2} E_2^-$, and
$E_n^+$ has measure zero for each odd integer $n$ except $-1$, and
$E_n^-$ has measure zero for each even integer $n$ except $0,2$.
Thus, $F$ must have the following form
\[ F = (E_{-2}^+ -4\pi) \cup E_0^- \cup ( E_{-1}^+ -2\pi) \cup
( E_1^- +2\pi ) \cup E_0^+ \cup (E_2^- + 4\pi). \]
 Using the fact that $E_1^- +2\pi = \frac{1}{2} E_{-2}^+$ and $E_{-1}^+
-2\pi = \frac{1}{2} E_2^-$, it is easy to verify that $(\sigma_E^F)^2 =
id_{\mathbb{R}}$.
\ \ \ \ \ \ \ \ \ \ \ \ \ \ \ \ \ \ \ \ \ \ \ \ \ \ \ \ \ \ \ \ \ \ \ \
\ \ \ \ \ \ \ \ \ \ \ \ \ \ \ \ \ \ \ \ \ \ \ \ \ \ \ \ \ \ \ \ \ \ \ \
\ \ \ \ \ \ \ \ \ \ \ \ \ \ \ \ \ \ \ \ \ \ \ \ \ \ \ \ \ \ \ \ \ \ \ \
$\square$
\\

\section{Some examples of wavelet sets}

Wavelet sets are useful as examples and counterexamples. Many
examples of them in the real line were given in [DL] exactly for
that purpose, for experimentation and testing hypotheses, and many
open questions still remain in that setting. The existence of
wavelet sets in the plane, and more generally in $\mathbb{R}^n$, was first
proven by Dai, Larson and Speegle [DLS1] in the summer of 1994
during a course at Texas A\&M taught by the second author using the
manuscript of [DL] as a text. The proof in [DLS1] was abstract and
covered dual congruence results for more general types of
dual-dynamical systems. But given the definitions, the proof was
constructive, implicitly giving an algorithm (a
$\emph{casting-outward}$ technique) for constructing examples of
wavelet sets for arbitrary expansive matrices and full rank
translation lattices on $\mathbb{R}^n$. However, this method produced only
wavelet sets that were unbounded and had $0$ as a limit point, and
were not very well described by diagrams or pictures. In other
words, they were wavelet sets but not "nice" wavelet sets. In 1995
concrete examples of dyadic wavelet sets in the plane, which were
bounded and bounded away from $0$, and could be nicely diagrammed,
were constructed by Soardi and Weiland [SW]. At about the same
time, Dai and Larson constructed two $\mathbb{R}^2$ dyadic wavelet sets
(denoted the "four-corners set" and the "wedding cake set"), for
inclusion in a final version of [DL] in response to a referee's
suggestion. In 1996-97, a number of authors constructed many other
concrete examples of wavelet sets in the plane, and in higher
dimensions. Included are the wavelet sets computed and
diagramed in the articles [BMM], [BL1], [BL2], [DLS2], [Za],
[C], [GH], [Gu2], [Han2], [S2]. The last three are the Ph.D. theses of
Gu, Han and Speegle, respectively. Importantly,
Baggett-Medina-Merrill [BMM], and Benedetto-Leon [BL2],
independently found two interesting completely different
constructive characterizations of all wavelet sets for expansive
matrices in $\mathbb{R}^n$. Open questions remain, especially questions
concerning the existence of wavelet sets with special properties,
and algorithms for constructing special classes of such sets.

 In this section, we will apply the techniques introduced in Lemma 5 to construct
an unbounded symmetric wavelet set and counterexamples to two
related questions on wavelet sets, in addition to new constructions
of some known wavelet sets.
\\
\\
{\bf Example 7.} Let $G_n = [\frac{2^n-1}{2^n}\pi,
\frac{2^{n+1}-1}{2^{n+1}}\pi) + (2^{n+2}-2)\pi,
n> 0$, then
\[ G_n \subseteq [(2^{n+2}-2)\pi, (2^{n+2}-1)\pi) = 2^{n+1}
[\pi + \frac{2^n-1}{2^n}\pi, \pi+\frac{2^{n+1}-1}{2^{n+1}}\pi). \]
Observe that $\dot{\bigcup}_{n>0} ( G_n - (2^{n+2} - 2)\pi ) =
[\frac{\pi}{2}, \pi)$ and
$\dot{\bigcup}_{n>0} \frac{1}{2^{n+1}} G_n \subset
\dot{\bigcup}_{n>0} [\pi + \frac{2^n-1}{2^n}\pi,
\pi+\frac{2^{n+1}-1}{2^{n+1}}\pi) = [\frac{3}{2}\pi, 2\pi)$.
Let $E_0
= [0,\frac{\pi}{2})$, let $F_0 = [\pi, 2\pi) \,\backslash\,
(\dot{\bigcup}_{n=1}^{\infty} \frac{1}{2^{n+1}} G_n)$. Notice that since
$F_0 \supseteq [\pi, \frac{3}{2}\pi)$,
\[ E_0 +2\pi \subseteq 2F_0 \mbox{\, \, \, and \, \, \, }
   \frac{1}{4}F_0 \subseteq E_0 , \]
and $( E_0+2\pi ) \cap \frac{1}{4} F_0 = \emptyset$.
Let $G_0$ be the set obtained by applying Lemma 5.
Then $\dot{\bigcup}_{n=0}^{\infty} G_n$ is both $2\pi$-translation
congruent to $[0,\pi)$ and $2$-dilation congruent to $[\pi, 2\pi)$. Hence,
by symmetry,  $- ( \dot{\bigcup}_{n=0}^{\infty} G_n) \,\dot{\cup}\, (
\dot{\bigcup}_{n=0}^{\infty} G_n )$ is an unbounded symmetric wavelet set.
\ \ \ \ \ \ \ \ \ \ \ \ \ \ \ \ \ \ \ \ \ \ \ \ \ \ \ \ \ \ \ \ \ \ \ \
\ \ \ \ \ \ \ \ \ \ \ \ \ \ \ \ \ \ \ \ \ \ \ \ \ \ \ \ \ \ \ \ \ \ \
\ \ \ \ \ \ \ \ \ \ \  \ \ \ \ \ \ \ \ \ \ \ \ \ \ \
$\square$
\\

We note that the first example of an unbounded symmetric wavelet set
was given in Proposition 2.14 of [FW]. The construction involved
some significant computations and explanation. Example 4.5 (xi) of
[DL] gave a different unbounded wavelet set whose construction
required little explanation, but it was not symmetric.

The following is a counterexample to a question posed by the second
author in [La1]:

{\it Let $E$ be a wavelet set, and suppose that $G \subseteq 2E\cup
E\cup \frac{1}{2} E$ is also a wavelet set. Is $(E, G)$ an
interpolation pair?}
\\
\\
{\bf Example 8.} Let $E=[-2\pi, -\pi)\cup [\pi, 2\pi)$, and $G\subseteq
2E \cup E \cup \frac{1}{2}E = [-4\pi, -\frac{\pi}{2}) \cup$
$[\frac{\pi}{2}, 4\pi)$ be a wavelet set. Suppose that $G_1 = G\cap
[3\pi, 4\pi)$ has positive measure,
then $G_1-2\pi \subseteq [\pi, 2\pi) \subseteq E$. $\forall s \in
G_1-2\pi$, $(\sigma_E^G)^2(s) = \sigma_E^G(s+2\pi) = 2\cdot
\sigma_E^G(\frac{1}{2}s+\pi) = 2\cdot (\frac{1}{2} s +\pi + 2k\pi) =
s + 2\pi + 4k\pi$, for some $k\in
\mathbb{Z}$. Thus, $(\sigma_E^G)^2 (s) \neq s$ for each $s\in
G_1 -2\pi$. Therefore, $(E,G)$ is not an interpolation pair if $G\cap
[3\pi, 4\pi)$ has positive measure. Since
\[ \frac{1}{2} [\pi, \frac{7}{4}\pi) \subseteq [0,\pi)
      \mbox{\, \, \, and \, \, \, }
   [0,\pi) +2\pi \subseteq 2[\pi, \frac{7}{4}\pi) , \]
and $\frac{1}{2} [\pi, \frac{7}{4}\pi) \,\cap\, ( [0,\pi) +2\pi )=
\emptyset$,
by Lemma 5, there exists $G_0 \subseteq [\frac{\pi}{2},
\frac{7}{8}\pi) \cup [2\pi, 3\pi)$, such that $G_0$ is
$2\pi$-translation congruent to $[0,\pi)$ and $2$-dilation congruent to
$[\pi,\frac{7}{4}\pi)$. Let $G = [-\pi, -\frac{\pi}{2}) \cup G_0
\cup [\frac{7}{2}\pi, 4\pi)$. Then $G$ is a wavelet set and $G \subseteq 2E\cup
E\cup \frac{1}{2} E$, but $(E, G)$ is not an interpolation pair.
\ \ \ \ \ \ \ \ \ \ \ \ \ \ \ \ \ \ \ \ \ \ \ \ \ \ \ \ \ \ \
\ \ \ \ \ \ \ \ \ \ \ \ \ \ \ \ \ \ \ \ \ \ \ \ \ \ \ \ \ \ \
$\square$
\\

The next example answers another question from [La1]:

{\it Let $E$ be a wavelet set, and suppose that $G \subseteq
(E-2\pi) \cup E \cup (E +2\pi) $ is also a wavelet set. Is $(E, G)$
an interpolation pair?}
\\
\\
{\bf Example 9.} Let $E = [-2\pi, -\pi) \cup [\pi, 2\pi)$ and $G
\subseteq (E-2\pi) \cup E \cup (E+2\pi) = [-4\pi, -3\pi) \cup
[-2\pi, 2\pi) \cup [3\pi, 4\pi)$ be a wavelet set. Suppose that
$G\cap [-\pi, \pi)$ is not a null set. Without loss of generality,
assume that $G_1 = G \cap [0,\pi)$ has positive measure. Then, for
each $s \in G_1 -2\pi \subset E$, $(\sigma_E^G)^2 (s) = \sigma_E^G
(s+2\pi) = 2^{-k} \cdot \sigma_E^G(2^k (s+2\pi))$, for some $k>0$
such that $2^k (s+2\pi) \in [\pi, 2\pi)$. Then $(\sigma_E^G)^2 (s) =
2^{-k} \cdot (2^k (s+2\pi) +2 n\pi) = s + 2\pi + 2^{-k} \cdot 2n\pi$
for some $n\in \{ 1, 0,-1\}$.  Thus, $(\sigma_E^G)^2 (s) \neq s$ for
each $s\in G_1 - 2\pi$. Therefore, if $G \cap [-\pi, \pi)$ is not a
null set, $(E, G)$ is not an interpolation pair. Based on the above
observation, let $G_0$ be the measurable set given by Lemma 5 and
the following two containment relations:
\[ [-2\pi, -\frac{3}{2}\pi) \subseteq [-2\pi, -\frac{11}{8}\pi) \cup
[-\frac{9}{8}\pi, -\pi), \]
\[ ( [-2\pi, -\frac{11}{8}\pi) \cup [-\frac{9}{8}\pi,
-\pi) ) -2\pi \subseteq 2 [-2\pi, -\frac{3}{2}\pi). \]
Then, $G := [\frac{5}{8}\pi, \frac{7}{8}\pi) \cup [\pi,
\frac{5}{4}\pi) \cup [\frac{7}{2}\pi, 4\pi) \cup [-\frac{3}{4}\pi,
-\frac{\pi}{2}) \cup G_0$ is a wavelet set, since
\begin{eqnarray*}
[0, 2\pi) & = &
 ([-2\pi, -\frac{11}{8}\pi) + 2\pi ) \,\dot{\cup}\, [\frac{5}{8}\pi,
\frac{7}{8}\pi) \,\dot{\cup}\, ( [-\frac{9}{8}\pi, -\pi) +
2\pi) \\
 &   & \,\dot{\cup}\, [\pi, \frac{5}{4}\pi) \,\dot{\cup}\, (
[-\frac{3}{4}\pi, -\frac{\pi}{2}) + 2\pi ) \,\dot{\cup}\, ([\frac{7}{2},
4\pi) -2\pi) ,
\end{eqnarray*}
and
\begin{eqnarray*}
[-2\pi, -\pi)  \,\cup\, [\frac{5}{8}\pi, \frac{5}{4}\pi)
& = &
[-2\pi, -\frac{3}{2}\pi) \,\dot{\cup}\, 2\cdot [-\frac{3}{4}\pi,
-\frac{\pi}{2}) \\
 & &  \,\dot{\cup}\, [\frac{5}{8}\pi,
\frac{7}{8}\pi)  \,\dot{\cup}\, \frac{1}{4} \cdot [\frac{7}{2}\pi,
4\pi)  \,\dot{\cup}\, [\pi, \frac{5}{4}\pi)
\end{eqnarray*}
is a $2$-dilation generator of a measurable partition of $\mathbb{R}$.
However, $G\cap [-\pi, \pi)$ is not a null set, hence $(\sigma_E^G)^2 \ne
id_{\mathbb{R}}$.
\ \ \ \ \ \ \ \ \ \ \ \ \ \ \ \ \ \ \ \ \ \ \ \ \ \ \ \ \ \ \ \ \ \ \ \
\ \ \ \ \ \ \ \ \ \ \ \ \ \ \ \ \ \ \ \ \ \ \ \ \ \ \ \ \ \ \ \ \ \ \ \
$\square$
\\

Lemma 5 can be useful in constructing special wavelet sets, and in
particular, parametric families of wavelet sets.
\\
\\
{\bf Example 10.} Let $l,m,n\in \mathbb{N}$ be given. Let
$E_+(m) = [0, 2^{-m+1}\pi)$, $F_+ = [\pi ,2\pi)$. Observe
that
\[  \frac{1}{2^m}F_+ \subseteq E_+(m)  \mbox{\, \, \, and \, \,
\, }
    E_+(m) + 2^l\pi \subseteq 2^l F_+ , \]
and $\frac{1}{2^m} F_+ \cap (E_+(m) +2^l \pi) = \emptyset$. Let
$G_+(l,m)$ be the measurable set obtained by applying Lemma 5.
Straightforward computation shows that
\[ G_+(l,m) = (\frac{2^l}{2^{l+m}-1}
\pi , 2^{-m+1}\pi) \cup [2^l\pi, 2^l\pi + \frac{2^l}{2^{l+m}-1} \pi ) .\]
Similarly, let $E_-(m) = [-2\pi+2^{-m+1}\pi, 0)$, $F_-(m) =
[-2\pi+2^{-m+1}\pi, -\pi+2^{-m}\pi)$. Observe that
\[  F_-(m) \subseteq E_-(m)  \mbox{\, \, \, and \, \, \, }
    E_-(m) - 2^{m+n}\pi+2^n\pi \subseteq 2^{m+n} F_-(m)  . \]
By applying Lemma 5 again, we get
\begin{eqnarray*}
G_-(m,n) & = & [ -2^{m+n}\pi + 2^n\pi -
\frac{2^{m+n}-2^n}{2^{m+n}-1}\pi, -2^{m+n}\pi+2^n\pi) \\
 & \cup &
[-2\pi+2^{-m+1}\pi, -\frac{2^{m+n}-2^n}{2^{m+n}-1}\pi ) .
\end{eqnarray*}
Thus, for fixed
$l,m,n$, $G_-(m,n) \cup G_+(l,m)$ is both $2\pi$-translation congruent
to $E_-(m) \cup E_+(m) = [-2\pi + 2^{-m+1}\pi, 2^{-m+1}\pi)$ and
$2$-dilation congruent to $F_-(m) \cup F_+$, hence is a wavelet
set.
\ \ \ \ \ \ \ \ \ \ \ \ \ \ \ \ \ \ \ \ \ \ \ \ \ \ \ \ \ \ \ \ \ \ \ \
\ \ \ \ \ \ \ \ \ \ \ \ \ \ \ \ \ \ \ \ \ \ \ \ \ \ \ \ \ \ \ \ \ \ \ \
\ \ \ \ \ \ \ \
$\square$
\\

Notice that $G_-(l,m,n) \cup G_+(l,m,n), l,m,n\ge 1$ is exactly
the family of wavelet sets $K_{l,m,n}, l,m,n\ge 1$ introduced by X. Fang
and X. Wang in [FW] (see [FW] Example (5)).

In the following examples, we will construct some known and new wavelet
sets in ${\mathbb{R}}^2$ with dilation matix $A$ given by $2\cdot
id_{{\mathbb{R}}^2}$. It is known (see [DLS1] and [SW]) that $E\subseteq
{\mathbb{R}}^2$ is a wavelet set for $2\cdot id_{{\mathbb{R}}^2}$ if and
only if $E$ is both $2\pi$-translation congruent to $[-\pi,\pi)\times
[-\pi, \pi)$ and $2$-dilation congruent to $[-2\pi,2 \pi)\times
[-2\pi, 2\pi) \,\backslash\, [-\pi,\pi)\times [-\pi, \pi)$.
Wavelet sets in higher dimensional spaces with
respect to arbitrary real expansive matrices can also be obtained in the
similar way.
\\
\\
{\bf Example 11.} Consider the first quadrant. Since we have the
following two relations:
\[ [0,\pi)\times [0,\pi) \,\backslash\, [0,\frac{\pi}{2}) \times [0,
\frac{\pi}{2}) \subseteq [0,\pi)\times [0,\pi), \]
\[ [0,\pi)\times [0,\pi) + (2\pi, 2\pi) \subseteq 4 \cdot  [0,\pi)\times
[0,\pi) \,\backslash\, [0,\frac{\pi}{2}) \times [0, \frac{\pi}{2}),  \]
and $[0,\pi)\times [0,\pi) \,\backslash\, [0,\frac{\pi}{2}) \times [0,
\frac{\pi}{2}) \,\cap\, ( [0,\pi)\times [0,\pi) + (2\pi, 2\pi) ) =
\emptyset$. Applying Lemma 5, we can construct a set $W_1$ which is both
$2\pi$-translation congruent to $[0,\pi)\times [0,\pi)$ and $2$-dilation
congruent to $[0,\pi)\times [0,\pi) \,\backslash\,
[0,\frac{\pi}{2}) \times [0, \frac{\pi}{2})$. Symmetrically, construct
sets $W_2, W_3, W_4$ in second, third and fourth quadrants, respectively.
Then $W_1 \cup W_2 \cup W_3 \cup W_4$ is a wavelet set. Straightforward
compuation shows that it is exactly the ``four corners set'' in [DLS2].
\ \ \ \ \ \ \ \ \ \ \ \ \ \ \ \ \ \ \ \ \ \ \ \ \ \ \ \ \ \ \ \ \
\ \ \ \ \ \ \ \ \ \ \ \ \ \ \ \ \ \ \ \ \ \ \ \ \ \ \ \ \ \ \ \ \
\ \ \ \ \ \ \ \ \ \ \ \ \ \ \ \ \
$\square$
\\
\\
{\bf Example 12.} Consider the right half plane. Since we have the
following two relations:
\[ [0,\pi)\times [-\pi,\pi) \,\backslash\, [0,\frac{\pi}{2}) \times
[-\frac{\pi}{2}, \frac{\pi}{2}) \subseteq [0,\pi)\times [-\pi,\pi), \]
\[ [0,\pi)\times [-\pi,\pi) + (2\pi, 0) \subseteq 4 \cdot
[0,\pi)\times [-\pi,\pi) \,\backslash\, [0,\frac{\pi}{2}) \times
[-\frac{\pi}{2}, \frac{\pi}{2}),  \] and $[0,\pi)\times [-\pi,\pi)
\,\backslash\, [0,\frac{\pi}{2}) \times [-\frac{\pi}{2},
\frac{\pi}{2}) \,\cap\, ( [0,\pi)\times [-\pi ,\pi) + (2\pi, 0) ) =
\emptyset$. By Lemma 5, we can construct a set $W_1$ which is both
$2\pi$-translation congruent to $[0,\pi)\times [-\pi,\pi)$ and
$2$-dilation congruent to $[0,\pi)\times [-\pi,\pi) \,\backslash\,
[0,\frac{\pi}{2}) \times [-\frac{\pi}{2}, \frac{\pi}{2})$.
Symmetrically, construct the set $W_2$ in the left half plane. Then
$W_1 \cup W_2$ is a wavelet set. Straightforward computation shows
that it is exactly the ``wedding cake set'' (Example 6.6.2 in [DL]
and also Figure 2 in[DLS2].) \ \ \ \ \ \ \ \ \ \ \ \ \ \ \ \ \ \ \ \
\ \ \ \ \ \ \ \ \ \ \ \ \ \ \ \ \ \ \ \ \ \ \ \ \ \ \ \ \ \ \ \ \ \
\ \ \ \ \ \ \ \ \ \ \ \ \ \ \ \ \ \ \ \ \ \ \ \ \ \ \ \ \ \ \ \ \ \ \ $\square$
\\
\\
{\bf Example 13.} Consider the left-top half plane (above the line
$y = x$ in the left half plane). Let $E_1 = \{ (x,y) | x\ge -\pi,
y\le \pi, y\ge x \}$, and $F_1 = \{ (x,y) | x\ge -\pi, y\le \pi,
y\ge x \} \,\,\backslash\,\, \{ (x,y) | x\ge -\frac{\pi}{2}, y\le
\frac{\pi}{2}, y\ge x \}$. Then we have the following two relations:
\[ F_1 \subseteq E_1, \ \ \ \ \ E_1 +(-2\pi, 2\pi) \subseteq 4\cdot
F_1, \] and $ F_1 \,\cap\,  ( E_1 +(-2\pi, 2\pi) ) = \emptyset$. Now
use Lemma 5 to construct a set $W_1$ which is both
$2\pi$-translation congruent to $E_1$ and $2$-dilation congruent to
$F_1$. Symmetrically, construct a corresponding set $W_2$ in the
right-bottom half plane. Then $W_1 \cup W_2$ is a wavelet set. The
diagram for this is given in Figure 1. \ \ \ \ \ \ \ \ \ \ \ \ \ \ \
\ \ \ \ \ \ \ \ \ \ \ \ \ \ \ \ \ \ \ \ \ \ \ \ \ \ \ \ \ \ \ \ \ \
$\square$
\\

\begin{center}
\quad\quad\epsfig{file=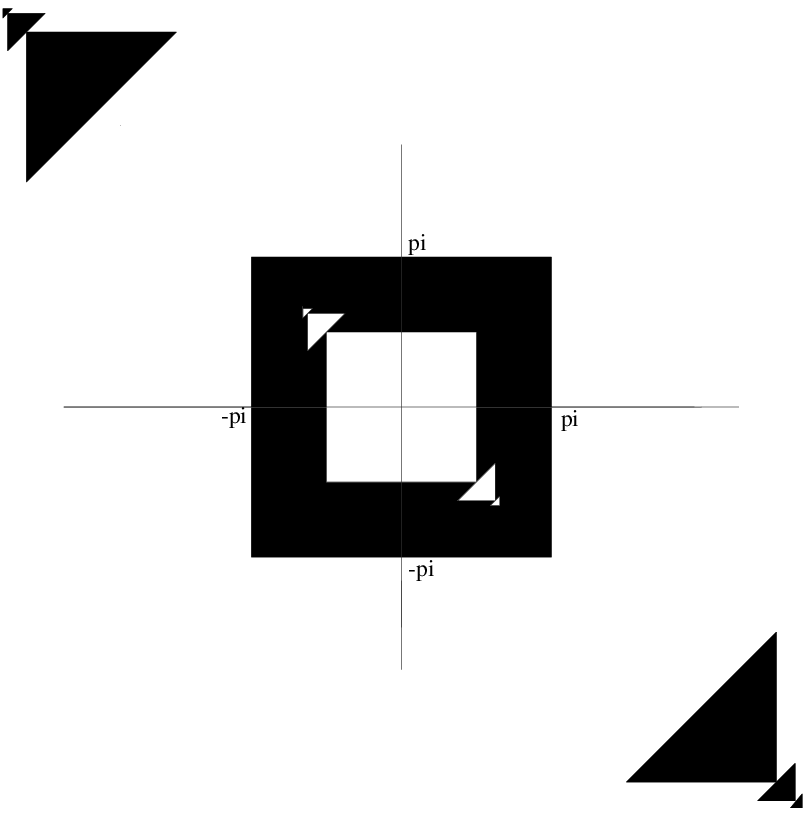,height=4.4in,width=3.6in,angle=0}
\vspace{-60pt}
\end{center}
\begin{center}
{\scshape Figure 1} Pine Tree
\end{center}

\vspace{15pt}

\hspace{-12pt}{\bf Example 14.} Consider the first quadrant. Let
$E_1 = [-\frac{\pi}{2}, 0) \times [-\frac{\pi}{2}, 0) \cup [-\pi,
-\frac{3}{4}\pi) \times [-\pi, -\frac{3}{4}\pi)$ and $F_1 =
[\frac{3}{2}\pi, 2\pi) \times [\frac{3}{2}\pi, 2\pi)$. Then we have
the following two relations:
\[ F_1 \subseteq E_1 + (2\pi, 2\pi), \ \ \ \ \ E_1 +(4\pi,
4\pi) \subseteq 2\cdot F_1, \]
and $ F_1 \,\cap\,  ( E_1 +(4\pi, 4\pi) ) = \emptyset$. By Lemma 5, we
can construct a set $W_1$ which is both $2\pi$-translation congruent to
$E_1$ and $2$-dilation congruent to $F_1$. Symmetrically, we can define
$E_2, E_3, E_4$ and $F_2, F_3, F_4$
, and  construct sets $W_2, W_3, W_4$ in the second, third and fourth
quadrants, respectively.
Let $B = [-\pi,\pi) \times [-\pi,\pi)
\,\backslash\, ( [-\frac{1}{2}\pi,\frac{1}{2}\pi) \times [-\frac{1}{2}\pi,\frac{1}{2}\pi)
\cup [\frac{3}{4}\pi,\pi) \times [\frac{3}{4}\pi,\pi)
\cup [-\pi,-\frac{3}{4}\pi) \times [\frac{3}{4}\pi,\pi)
\cup [\frac{3}{4}\pi,\pi) \times [-\pi,-\frac{3}{4}\pi)
\cup [[-\pi,-\frac{3}{4}\pi) \times [-\pi,-\frac{3}{4}\pi) )$.
Then since $B \cup E_1 \cup E_2 \cup E_3 \cup
E_4 = [-\pi, \pi)\times[-\pi,\pi)$ and $2B \cup F_1 \cup F_2 \cup F_3 \cup
F_4 = [-2\pi, 2\pi)\times[-2\pi,2\pi) \,\backslash\, [-\pi,
\pi)\times[-\pi,\pi)$, $W_1 \cup W_2 \cup W_3 \cup W_4 \cup B$
is a wavelet set. Computation shows that it is one of the
wavelet sets introduced in [SW]. The diagram for this is given in Figure
2.
\ \ \ \ \ \ \ \ \ \ \ \ \ \ \ \ \ \ \ \ \ \ \ \ \ \ \ \ \ \ \ \ \ \
\ \ \ \ \ \ \ \ \ \ \ \ \ \ \ \ \ \ \ \ \ \ \ \ \ \ \
$\square$

\vspace{15pt}

\begin{center}
\quad\quad\epsfig{file=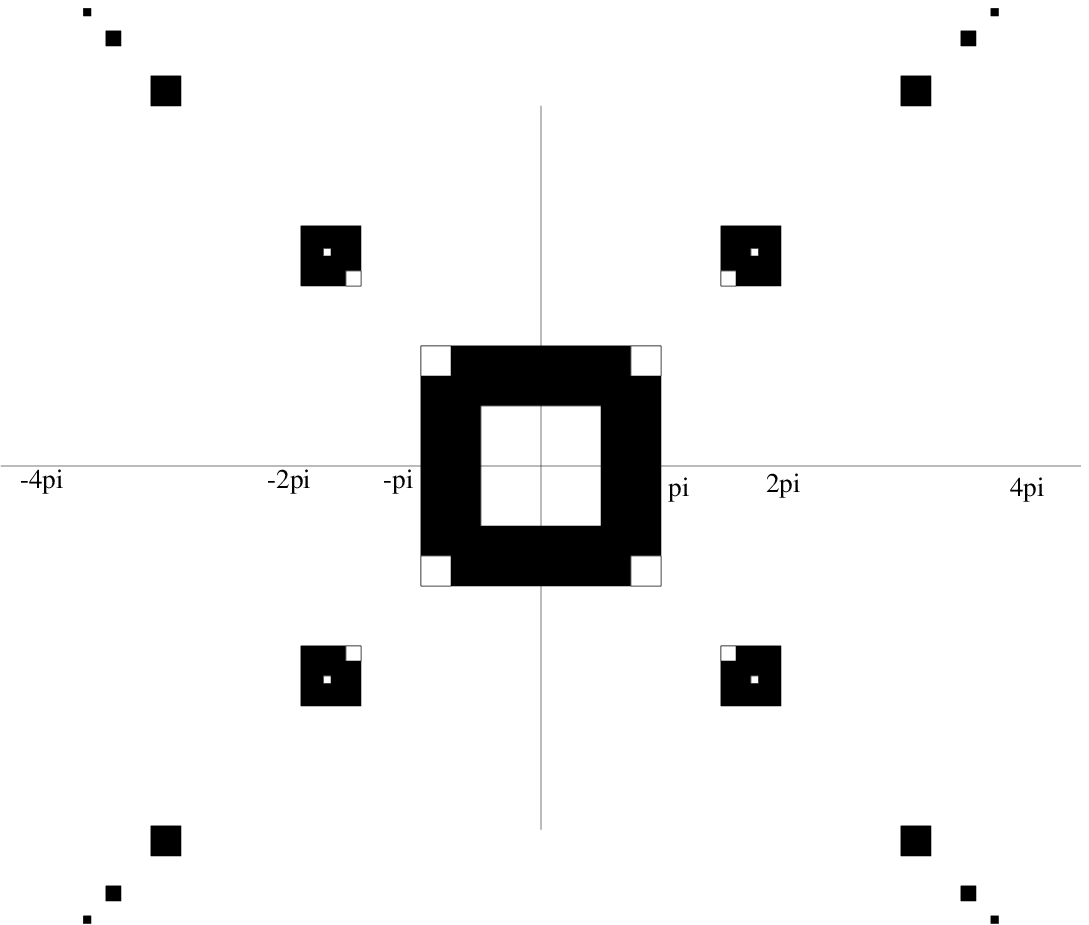,height=3.6in,width=4in,angle=0}
\end{center}

\vspace{10pt}
\begin{center}
{\scshape Figure 2} A Wavelet Set in [SW]
\end{center}

\vspace{10pt}

Remarks. (1) The idea of Lemma 5 was used in
constructing a covering of $\mathbb{R}$ by symmetric wavelet sets,
and this was a key to subsequent work of the first author with
Rzeszotnik in [RZ].

(2) Further Directions: Interpolation maps, and interpolation pairs
and more general interpolation families of wavelet sets, make sense
and have been studied for matrix dilations in ${\mathbb{R}}^n$.
Congruence domains also make sense for matrix dilations.  Lemma 5 in
this article was stated and proved for ${\mathbb{R}}^n$, and was
applied to solve Question B in ${\mathbb{R}}^1$, and also used to
study examples of dyadic (i.e. for matrix dilation $2I$) wavelet
sets in the plane. The results in this paper suggest some directions
for further research. In particular, does Theorem 1 extend to matrix
dilations in ${\mathbb{R}}^n$, especially for the cases where
interpolation pairs are known to exist? It might be useful to try to
extend Lemma 5 and also Theorem 3 to matrix dilations.

\renewcommand{\baselinestretch}{1}

\bibliographystyle{alpha}

\end{document}